\documentclass[12pt,a4paper,reqno]{amsart}
\usepackage[english]{babel}
\usepackage{amssymb,latexsym,amsfonts,amsthm,upref,amsmath}
\usepackage[foot]{amsaddr}
\usepackage[margin=1in]{geometry}
\usepackage[dvipsnames,x11names]{xcolor}
\definecolor{myblue}{HTML}{003399}
\usepackage{hyperref}
\hypersetup{colorlinks,citecolor=myblue,filecolor=black,linkcolor=myblue,urlcolor=myblue}
\usepackage{enumerate} 
\usepackage{tikz}
\usepackage{float}
\usepackage{cleveref}
\usepackage{mathtools}
\usepackage{cite}

\usepackage{filecontents}
\makeatletter
\newcommand{\leqnomode}{\tagsleft@true}
\newcommand{\reqnomode}{\tagsleft@false}
\newcommand{\cev}[1]{\reflectbox{\ensuremath{\vec{\reflectbox{\ensuremath{#1}}}}}}
\makeatother

\newtheorem*{corintro*}{Corollary}
\newtheorem*{thm*}{Theorem}
\newtheorem*{lem*}{Lemma}
\newtheoremstyle{prim}{}{}{\normalfont}{}{\bfseries}{.}{ }{}
\theoremstyle{prim}

\newtheoremstyle{stil}{}{}{\slshape}{}{\bfseries}{.}{ }{}
\theoremstyle{stil}
\newtheorem{thm}{Theorem}[section]
\newtheoremstyle{defi}{}{}{}{}{\bfseries}{.}{ }{}
\theoremstyle{defi}

\theoremstyle{defi}
\newtheorem{rem}[thm]{Remark}
\theoremstyle{stil}
\newtheorem{pro}[thm]{Proposition}
\theoremstyle{stil}
\newtheorem{lem}[thm]{Lemma}
\theoremstyle{stil}
\newtheorem{kor}[thm]{Corollary}
\newenvironment{prf}{\noindent \textit{Proof.}}{\null\hfill$\qed$\hskip
2mm\vskip 2mm}

\newcommand{\sgn}{ {\rm sgn}\ts}

\newcommand{\vac}{\mathop{\mathrm{\boldsymbol{1}}}}
\newcommand{\ndo}{\mathop{\mathrm{End}}}
\newcommand{\om}{\mathop{\mathrm{Hom}}}

\newcommand{\grrr}{\mathop{\mathrm{gr}}}

\newcommand{\Sym}{\mathfrak S}

\newcommand{\DY}{ {\rm DY}}
\newcommand{\Y}{ {\rm Y}}

\newcommand{\gl}{\mathfrak{gl}}
\newcommand{\sll}{\mathfrak{sl}}

\newcommand{\cdotrl}{\mathop{\hspace{-2pt}\underset{\text{RL}}{\cdot}\hspace{-2pt}}}
\newcommand{\cdotlr}{\mathop{\hspace{-2pt}\underset{\text{LR}}{\cdot}\hspace{-2pt}}}

\newcommand{\wtld}{\widetilde}
\newcommand{\wht}{\widehat}
\newcommand{\wvr}{\overline}
\newcommand{\wndr}{\underline}

\newcommand{\ot}{\otimes}

\newcommand{\ts}{\,}
\newcommand{\tss}{\hspace{1pt}}

\newcommand{\CC}{\mathbb{C}\tss}

\newcommand{\ZZ}{\mathbb{Z}\tss}
\newcommand{\TT}{\mathbb{T}}
\newcommand{\BB}{\mathbb{A}}
\newcommand{\sss}{\mathbb{A}}
\newcommand{\mmm}{\mathbb{M}}

\newcommand{\Ac}{ {\rm A}}

\newcommand{\Sc}{\mathcal{S}}

\newcommand{\Ec}{\mathcal{E}}

\newcommand{\Uc}{\mathcal{U}}
\newcommand{\Vc}{\mathcal{V}}
\newcommand{\Wc}{\mathcal{M}}

\newcommand{\g}{\mathfrak{g}}
\newcommand{\oo}{\mathfrak{o}}
\newcommand{\spp}{\mathfrak{sp}}

\newcommand{\z}{\mathfrak{z}}

\newcommand{\sdet}{ {\rm sdet}\ts}
\newcommand{\qdet}{ {\rm qdet}\ts}
\newcommand{\tr}{ {\rm tr}}
\newcommand{\fand}{\quad\text{and}\quad}
\newcommand{\Fand}{\qquad\text{and}\qquad}

\newcommand{\non}{\nonumber}
\newcommand{\beq}{\begin{equation}}
\newcommand{\eeq}{\end{equation}}
\newcommand{\ben}{\begin{equation*}}
\newcommand{\een}{\end{equation*}}

\makeatletter
\def\smalloverbrace#1{\mathop{\vbox{\m@th\ialign{##\crcr\noalign{\kern3\p@}%
  \tiny\downbracefill\crcr\noalign{\kern3\p@\nointerlineskip}%
  $\hfil\displaystyle{#1}\hfil$\crcr}}}\limits}
\makeatother

\makeatletter
\def\smallunderbrace#1{\mathop{\vtop{\m@th\ialign{##\crcr
   $\hfil\displaystyle{#1}\hfil$\crcr
   \noalign{\kern3\p@\nointerlineskip}%
   \tiny\upbracefill\crcr\noalign{\kern3\p@}}}}\limits}
\makeatother

\begin{document}

\title[A note on constructing quasi modules for quantum vertex algebras  from twisted Yangians]{A note on constructing quasi modules for quantum vertex algebras  from twisted Yangians}

\author{Slaven Ko\v{z}i\'{c}} 
\address[S. Ko\v{z}i\'{c}]{Department of Mathematics, Faculty of Science, University of Zagreb, Bijeni\v{c}ka cesta 30, 10\,000 Zagreb, Croatia}
\email{kslaven@math.hr}

\author{Marina Serti\'{c}} 
\address[M. Serti\'{c}]{Department of Mathematics, Faculty of Economics \& Business, University of Zagreb, Trg J. F. Kennedy 6, 10\,000 Zagreb, Croatia}
\email{msertic1@net.efzg.hr}



\begin{abstract}
In this note, we consider the twisted Yangians  $\text{Y}(\mathfrak{g}_N)$ associated with the orthogonal and symplectic Lie algebras $\mathfrak{g}_N=\mathfrak{o}_N,\mathfrak{sp}_N$. First, we introduce a certain subalgebra  $\text{A}_c(\mathfrak{g}_N)$ of the  double Yangian for $\mathfrak{gl}_N$ at the level $c\in\mathbb{C}$, which contains   the  centrally extended $\text{Y}(\mathfrak{g}_N)$ at the level $c$ as well as its vacuum module  $\mathcal{M}_c(\mathfrak{g}_N)$. Next, we employ  its structure to construct examples of quasi modules for the quantum affine vertex algebra $\mathcal{V}_c(\mathfrak{gl}_N)$ associated  with the Yang $R$-matrix.
Finally, we use the  description of the center of $\mathcal{V}_c(\mathfrak{gl}_N)$ to obtain  explicit formulae for families of central elements for a certain   completion of $\text{A}_c(\mathfrak{g}_N)$ and invariants of $\mathcal{M}_c(\mathfrak{g}_N)$. 
\end{abstract}

\maketitle

\section{Introduction}\label{sec0}
The twisted Yangians 
are certain subalgebras of the Yangian for $\mathfrak{gl}_N$ which are
associated with the orthogonal and symplectic Lie algebras. 
They were introduced and studied by G. Olshanski in  \cite{O}. Later on, their properties were further investigated by A. Molev, M. Nazarov and G. Olshanski  \cite{MNO}. In particular, the Sklyanin determinant,  
which is a twisted analogue of the quantum determinant for the Yangian for $\mathfrak{gl}_N$,
was defined and studied in \cite{O,MNO}; see also the papers \cite{M1,M2,M3,JZ}. It was named  after E. K. Sklyanin in recognition
of his  work \cite{S}, where the new type of   determinant  for a certain class of reflection algebras was introduced.
For more information on twisted Yangians and  the Sklyanin determinant see the book by A. Molev \cite{M}.

The goal of this paper is to establish a connection between  the quantum vertex algebra theory and the twisted Yangians for  $\mathfrak{g}_N=\mathfrak{o}_N,\mathfrak{sp}_N$. In contrast with the Etingof--Kazhdan construction, which can be used to associate quantum vertex algebras with the double Yangians of   classical types \cite{BJK, EK}, the vertex operators coming from the twisted Yangians  no longer possess the $\Sc$-locality property, a quantum analogue of the locality of vertex operators. Thus, we were only able to employ the structure of twisted Yangians to construct families of quasi modules for  Etingof--Kazhdan's quantum   vertex algebra $\mathcal{V}_c(\mathfrak{gl}_N)$ \cite{EK} associated  with the Yang $R$-matrix. We should say that the notion of quasi module, which was introduced by H.-S. Li \cite{Li},    presents a generalization  of the  vertex algebra module.

Our construction   goes in  parallel with the construction of quasi $\mathcal{V}_c(\mathfrak{gl}_N)$-modules associated with
reflection algebras  \cite{K19}, which is to be expected
given the numerous similarities    between the properties of these two classes of algebras. We start by introducing a certain subalgebra  $\text{A}_c(\mathfrak{g}_N)$ of the double Yangian for $\mathfrak{gl}_N$ at the level $c\in\mathbb{C}$, which contains both  the suitably defined central extension  of the twisted Yangian $\text{Y}(\mathfrak{g}_N)$ at the level $c$ as well as its vacuum module  $\mathcal{M}_c(\mathfrak{g}_N)$. The action of $\text{A}_c(\mathfrak{g}_N)$ is then used to equip the $\CC[[h]]$-module of $\mathcal{M}_c(\mathfrak{g}_N)$  with the structure of quasi $\mathcal{V}_{2c}(\mathfrak{gl}_N)$-module so that, in particular, the action of the twisted Yangian   resembles the annihilation operators. In addition, we show that a certain wide class of $\text{A}_c(\mathfrak{g}_N)$-modules, the so-called restricted modules, is naturally equipped with the structure of quasi $\mathcal{V}_{2c}(\mathfrak{gl}_N)$-module.

Finally, as an application, we combine the corresponding quasi module maps with the   description of the center of  $\mathcal{V}_{2c}(\mathfrak{gl}_N)$, as given by N. Jing, A. Molev, F. Yang and the first author \cite{JKMY}, to obtain explicit formulae for families of central elements of the  completed algebra $\wtld{\text{A}}_c(\mathfrak{g}_N)$ and invariants of the quasi module $\mathcal{M}_c(\mathfrak{g}_N)$. At the critical level $2c=-N$,  we recover large families of central elements and invariants parametrized by   certain Young diagrams, while at the noncritical level $2c\neq -N$, we obtain a vertex operator-theoretic interpretation of the   Sklyanin determinant.

\section{Preliminaries}\label{sec1}
\numberwithin{equation}{section}
In this section, we  recall the    double Yangian for $ \gl_N $ and   Etingof--Kazhdan's construction  of the quantum affine vertex algebra associated with the Yang $R$-matrix.

\subsection{Double Yangian for \texorpdfstring{$\gl_N$}{glN}}\label{sec11}
Let $N\geqslant 2$ be an   integer and   $h$   a formal parameter.
 The {\em Yang $R$-matrix} over the commutative ring $\CC[[h]]$ is defined by
\beq\label{yang}
R(u)=1-\textstyle\frac{h}{u}P .
\eeq
Here $1$ denotes the identity and $P$  the permutation operator on $\CC^N\ot \CC^N  $, i.e. we have
\beq\label{perm}
1=\sum_{i,j=1}^N e_{ii}\ot e_{jj} \fand P=\sum_{i,j=1}^N e_{ij}\ot e_{ji},
\eeq
where $e_{ij}$ are matrix units. 
Note that 
$P^{t_1}=P^{t_2}$,
where $t_n$ is the transposition $e_{rs}\mapsto e_{sr}$   applied on the $n$-th tensor factor.
The
$R$-matrix \eqref{yang} satisfies the {\em Yang--Baxter equation}
\beq\label{ybe}
R_{12}(u)\ts R_{13}(u+v)\ts R_{23}(v)
=R_{23}(v)\ts R_{13}(u+v)\ts R_{12}(u).
\eeq
Let $g(u)$  be
the unique   series in $1+u^{-1}\CC[[u^{-1}]]$
such that
$g(u+N)=g(u)(1-u^{-2})$. 
Then the $R$-matrix $\wvr{R}(u)=g(u/h)R(u)$ satisfies the {\em crossing symmetry relation},
\beq\label{csym}
\left(\wvr{R}(u)^{-1}\right)^{t}\ts \wvr{R}(u+hN)^{t}=1,
\eeq
where $\wvr{R}(u)^{t}$ denotes the transposed $R$-matrix $\wvr{R}(u)^{t_1}=\wvr{R}(u)^{t_2}$,  
and the {\em unitarity relation},
\beq\label{uni}
\wvr{R}(u)\ts \wvr{R}(-u)=1,
\eeq
see, e.g., \cite[Sect. 2]{JKMY} for more details.

\allowdisplaybreaks

The {\em double Yangian} $\DY(\gl_N)$ for $\gl_N$  is defined as the associative algebra over $\CC[[h]]$ generated by the   central element $C$ and the elements
$t_{ij}^{(\pm r)}$, where $i,j=1,\ldots ,N$ and $r=1,2,\ldots ;$  see, e.g., \cite{I}.   Its defining relations are given by
\begin{align}
R\big(u-v \big)\ts T_1(u)\ts T_2(v)&=T_2(v)\ts T_1(u)\ts
R\big(u-v \big),\label{RTT2}\\
R\big(u-v \big)\ts T^+_1(u)\ts T^+_2(v)&=T^+_2(v)\ts T^+_1(u)\ts
R\big(u-v \big),\label{RTT1}\\
\wvr{R}\big(u-v+hC/2\big)\ts T_1(u)\ts T^+_2(v)&=T^+_2(v)\ts T_1(u)\ts
\wvr{R}\big(u-v-hC/2\big).\label{RTT3}
\end{align}
The matrices $T(u)$ and $T^+(u)$   are defined by
\begin{align*}
T(u)=\sum_{i,j=1}^N e_{ij}\ot t_{ij}(u)\Fand T^{+}(u)=\sum_{i,j=1}^N e_{ij}\ot t_{ij}^{+}(u),
\end{align*}
while its entries, the series $t_{ij}(u) $ and $t_{ij}^+ (u) $, are given by
$$t_{ij}(u)=\delta_{ij}+h\sum_{r=1}^{\infty} t_{ij}^{(r)}u^{-r}\Fand t_{ij}^+ (u)=\delta_{ij}-h\sum_{r=1}^{\infty} t_{ij}^{(-r)}u^{r-1}.$$
Throughout the paper we use the  subscript to indicate a copy of the matrix in the tensor product algebra
$ (\ndo\CC^N)^{\ot m}\ot\DY(\gl_N)$,
so that, for example,  we have
$$T_k(u)=\sum_{i,j=1}^{N} 1^{\otimes (k-1)} \ot e_{ij} \ot 1^{\ot (m-k)} \ot t_{ij}(u).$$
In particular, we have  $k=1,2$ and $m=2$  in the defining relations \eqref{RTT2}--\eqref{RTT3}.

Recall that the {\em $h$-adic topology} on an arbitrary $\CC[[h]]$-module $V$ is the topology generated by the  basis $v+h^n V$, $v\in V$, $n\in\mathbb{Z}_{\geqslant 0}$. We shall often write $V_h$ to indicate that the $\CC[[h]]$-module $V$ is $h$-adically completed. For example, if $W$ is a $\CC[[h]]$-module, then $V=W[z^{-1}]_h$ (resp. $V=W((z))_h$) denotes the $\CC[[h]]$-module of all power series $\sum_{r\in\ZZ} a_r z^r$ in $V$ such that $a_r\to 0$ when $r\to -\infty $ with respect to the $h$-adic topology.

The {\em Yangian} $\Y(\gl_N)$ (resp. the {\em dual Yangian } $\Y^+ (\gl_N)$) is the subalgebra of  $\DY(\gl_N)$ generated by the elements $t_{ij}^{(r)}$ (resp. $t_{ij}^{(-r)}$), where $i,j=1,\ldots,N$, $r=1,2,\ldots .$ 
For any  $c\in\CC$ we denote by $\DY_c (\gl_N)$ the {\em double Yangian at the level $c$}, i.e. the quotient of the algebra $\DY(\gl_N)$ by the ideal generated by $C-c$.
The {\em vacuum module $\Vc_c(\gl_N)$ at the level $c$} over the double Yangian is defined as the $h$-adic completion of the quotient of the algebra $\DY_c (\gl_N)$ by its left ideal generated by the elements $t_{ij}^{(r)}$, $i,j=1,\ldots ,N$, $ r=1,2,\ldots.$  
The Poincar\'e--Birkhoff--Witt theorem for the double Yangian, see  \cite[Thm. 2.2]{JKMY} or \cite[Thm. 15.3]{Naz}, implies that the  vacuum module   is isomorphic, as a $\CC[[h]]$-module, with the $h$-adically completed  dual Yangian.

\subsection{Quantum affine vertex algebra}

Let $n$ and $m$ be positive integers. For the families of variables $u=(u_1,\ldots,u_n)$ and $v=(v_{1},\ldots, v_{m})$ and the variable $z$ define
the formal power series with coefficients in  
$(\ndo\CC^N)^{\ot n} \ot (\ndo\CC^N)^{\ot m} $
by
\begin{align}
& \wvr{R}_{nm}^{12}(u|v|z)=\prod_{i=1,\ldots,n}^{\longrightarrow} \prod_{j=n+1,\ldots,n+m}^{\longleftarrow}\wvr{R}_{ij}(z+u_i-v_{j-n}),\label{bnm1}\\
&\wndr{\wvr{R}}_{nm}^{12}(u|v|z)=\prod_{i=1,\ldots,n}^{\longrightarrow} \prod_{j=n+1,\ldots,n+m}^{\longrightarrow}\wvr{R}_{ij}(z+u_i+v_{j-n}),\label{bnm2}
\end{align}
where the arrows indicate the order of the factors.
Moreover,  we  use the superscript $t$ to denote the product of transposed $R$-matrices, e.g.,
\beq\label{bnm3}
\wndr{\wvr{R}}_{nm}^{t\,12}(u|v|z)=\prod_{i=1,\ldots,n}^{\longrightarrow} \prod_{j=n+1,\ldots,n+m}^{\longrightarrow}\wvr{R}_{ij}^{t}(z+u_i+v_{j-n}).
\eeq
The series   $R_{nm}^{12}(u|v|z)$ and $\wndr{R}_{nm}^{12}(u|v|z)$  which employ the original (non-normalized) Yang $R$-matrix \eqref{yang} are introduced analogously.
In \eqref{bnm1}--\eqref{bnm3}, as well as in the rest of the paper, we use the common expansion convention, where the expressions of the form $(a_1 z_1+\ldots+ a_n z_n )^k$, with $a_i \in\CC$, $a_i\neq 0$ and $k<0$,   are expanded in negative powers of the variable which appears  on the left. Hence, for example, we have
$$(z_1 -z_2)^{-1}=\frac{1}{z_1}\sum_{l\geqslant 0}\left(\frac{z_2 }{z_1 }\right)^l \neq
\frac{1}{z_2}\sum_{l\geqslant 0}\left(\frac{z_1 }{z_2 }\right)^l=
 (-z_2+z_1)^{-1}.$$

Consider the following operators on $(\ndo\CC^N)^{\ot n} \ot\Vc_c(\gl_N)$:
\begin{align}
T_{[n]}^+(u|z)=T_1^+(z+u_1)\ldots T_n^+(z+u_n)\,\text{ and }\, T_{[n]}(u|z)=T_1(z+u_1)\ldots T_n(z+u_n).\label{teovi}
\end{align}
Using defining relations \eqref{RTT2}--\eqref{RTT3},  one can   verify the   equalities  from \cite[Subsect. 2.1.2]{EK},
\begin{align}
&R_{nm}^{12}(u|v|z-w)T_{[n]}^{+13}(u|z)T_{[m]}^{+23}(v|w)
=T_{[m]}^{+23}(v|w)T_{[n]}^{+13}(u|z)R_{nm}^{12}(u|v|z-w),\label{rtt1}\\ 
&R_{nm}^{12}(u|v|z-w)T_{[n]}^{13}(u|z)T_{[m]}^{23}(v|w)
=T_{[m]}^{23}(v|w)T_{[n]}^{13}(u|z)R_{nm}^{12}(u|v|z-w),\label{rtt2}\\ 
&\overline{R}_{nm}^{\ts 12}(u|v|z-w+h\tss c/2)T_{[n]}^{13}(u|z)T_{[m]}^{+23}(v|w) =T_{[m]}^{+23}(v|w)T_{[n]}^{13}(u|z)\overline{R}_{nm}^{\ts 12}(u|v|z-w-h\tss c/2),
\label{rtt3}
\end{align}
where the superscripts $1,2,3$ indicate the tensor product factors as follows:
\beq\label{fctr}
\smalloverbrace{(\ndo\CC^N)^{\ot n}}^{1} \ot \smalloverbrace{ (\ndo\CC^N)^{\ot m}}^{2}\ot  \smalloverbrace{\Vc_c(\gl_N)}^{3}.
\eeq 

 In      \eqref{bnm1}--\eqref{bnm3} and \eqref{teovi}, we sometimes omit the variable $z$ and then write, e.g.,
\beq\label{omt}
T_{[n]}^+(u)=T_1^+(u_1)\ldots T_n^+(u_n)\fand R_{nm}^{12}(u|v)=\prod_{i=1,\ldots,n}^{\longrightarrow} \prod_{j=n+1,\ldots,n+m}^{\longleftarrow}R_{ij}(u_i-v_{j-n}).
\eeq
For example, the identity in \eqref{rtt1}  with omitted variables   $z$ and $w$ takes  the form
$$
R_{nm}^{12}(u|v )T_{[n]}^{+13}(u )T_{[m]}^{+23}(v )
=T_{[m]}^{+23}(v )T_{[n]}^{+13}(u )R_{nm}^{12}(u|v ).
$$
Finally, we  recall   Etingof--Kazhdan's construction \cite[Thm. 2.3]{EK}:
\begin{thm}\label{EK:qva}
For any $c\in \CC$
there exists a unique  structure of quantum vertex algebra
on the $\CC[[h]]$-module $\Vc_c(\gl_N)$ such that the vacuum vector is
the unit $\vac\in \Vc_c(\gl_N)$ and the vertex operator map is given by
\beq\label{qva1}
Y\big(T_{[n]}^{+}(u )\vac,z\big)=T_{[n]}^{+}(u|z)\ts T_{[n]} (u|z+h\tss c/2)^{-1}.
\eeq
\end{thm}

 \section{Quasi modules for the quantum affine vertex algebra}\label{sec2} 
 In this section, we employ the   twisted Yangians  associated with the 
orthogonal and symplectic Lie algebras $\mathfrak{g}_N=\mathfrak{o}_N,\mathfrak{sp}_N$  to introduce certain  subalgebras $\Ac_c(\mathfrak{g}_N)$  of the double Yangian $\DY_c (\gl_N)$. Using these subalgebras we construct a   family of quasi $\Vc_c(\gl_N)$-modules and investigate  their connection  with the representation theory of  $\Ac_c(\mathfrak{g}_N)$. 

 \subsection{Twisted Yangians} 
We follow \cite[Ch. 2]{M} to define (slightly modified)  twisted Yangians. In addition, we introduce another family of related algebras which resemble  dual Yangians.
Let $G=(g_{ij})_{i,j=1}^N$ be a nonsingular complex matrix satisfying $G^t =\pm G$. If $G$ is symmetric (resp.  antisymmetric), denote by $\g_N$ the orthogonal Lie algebra $\mathfrak{o}_N$ (resp. the symplectic Lie algebra $\mathfrak{sp}_N$), where in the antisymmetric case   $N$ is even. 
Note that $\g_N$ coincides with the fixed point subalgebra of the   automorphism $\sigma   $  of $\mathfrak{gl}_N$ given by
\beq\label{sigma}
\sigma(A)=-G^{-1}  A^t \ts G\quad\text{for all}\quad A\in \mathfrak{gl}_N.
\eeq

Consider the series 
\beq\label{Sseries}
S^{+}(u)=T^{+}(u) \ts G \ts T^{+t}(-u)\fand S(u)=T(u+hc) \ts G \ts T^{t}(-u).
\eeq
They can be written in the form
$$
 S^{+}(u)=\sum_{i,j=1}^N e_{ij}\ot s_{ij}^{+}(u) \fand 
 S(u)=\sum_{i,j=1}^N e_{ij}\ot s_{ij}(u),$$
where the matrix entries are given by
\begin{align*}
 s_{ij}^{+}(u)=g_{ij} -h\sum_{r=1}^\infty s_{ij}^{(-r)}u^{r-1}\fand 
  s_{ij}(u)=g_{ij} +h\sum_{r=1}^\infty s_{ij}^{(r)}u^{-r} 
\end{align*}
for some  elements $s_{ij}^{(-r)}\in\Y^+(\gl_N)$ and $s_{ij}^{(r)}\in\Y(\gl_N)$.
The series \eqref{Sseries} satisfy
\begin{align}
&S^{+t}(-u)=\pm S^+(u) +\frac{h}{2u}\left(S^+(u)-S^+(-u)\right),\label{srel1}\\
&S^{t}(-u-hc)=\pm S(u) +\frac{h}{2u+hc} \left(S(u)-S(-u-hc)\right),\label{srel2}
\end{align}
where the plus (resp. minus) sign corresponds to the symmetric case  $\g_N =\oo_N$ (resp.  
 antisymmetric case $\g_N=\spp_N$). 
Indeed, both identities are easily verified using    \eqref{RTT2} and \eqref{RTT1}.
Furthermore,  the relations  \eqref{RTT2}--\eqref{RTT3}, along with   the equality
\beq\label{usflg}
R(u-v)\ts G_1\ts  R^t (-u-v)\ts  G_2 =G_2\ts  R^t (-u-v)\ts  G_1\ts  R(u-v) 
\eeq
from \cite[Lemma 2.4.1]{M},  imply the identities
\begin{align}
&R(u-v)\ts S_1^+(u)\ts  R^t(-u-v)\ts  S_2^+(v)=S_2^+(v)\ts R^t (-u-v)\ts S_1^+(u)\ts R(u-v)\label{RSRS1},\\
&R(u-v)\ts S_1(u)\ts  R^t (-u-v-hc)\ts  S_2(v)=S_2(v)\ts R^t(-u-v-hc)\ts S_1(u)\ts R(u-v)\label{RSRS2},\\
&\wvr{R}(u-v+3hc/2)\ts S_1(u)\ts  \wvr{R}^t(-u-v+hc/2) \ts S_2^+(v)\non\\
&\qquad
\qquad\qquad\qquad\qquad\ts\quad\ts
=S_2^+(v)\ts \wvr{R}^t(-u-v-3hc/2)\ts S_1(u)\ts \wvr{R}(u-v-hc/2).\label{RSRS3}
\end{align}

Let $\Ac_c  (\g_N)$ be  the subalgebra of the double Yangian  $\DY_c(\gl_N)$ at the level $c$   generated by the elements $s_{ij}^{(\pm r)}$  with $i,j=1,\ldots,N$ and $r=1,2,\ldots .$
Next,  let $\Y^+ (\g_N)$ (resp. $\Y_c  (\g_N)$) be the subalgebra  of the dual Yangian $\Y^+ (\mathfrak{gl}_N)$ (resp. the Yangian $\Y  (\mathfrak{gl}_N)$) generated by the elements $s_{ij}^{(-r)}$  (resp. $s_{ij}^{(r)}$) with $r\geqslant 1$. 
Finally, let $\Wc'_c (\g_N)$ be the   $\Ac_c(\g_N)$-submodule of $\Vc_c(\gl_N)$ generated by the vacuum vector $\vac $. 
Note that  $\Wc'_c (\g_N)$  coincides     with   $\Y^+(\g_N)$  as a $\CC[[h]]$-module. 
We define the {\em vacuum module} $\Wc_c (\g_N)$ for the algebra  $\Ac_c  (\g_N)$ 
as the $h$-adic completion of $\left\{v\in\Vc_c(\gl_N)   : h^n v\in  \Wc'_c (\g_N)\text{ for some }n\geqslant 0 \right\}$. By the definition,  the induced topology on $\Wc_c (\g_N)$  from  $\Vc_c(\gl_N)$ coincides with the $h$-adic topology of $\Wc_c (\g_N)$; cf. \cite[Lemma 3.5]{Li}. Moreover, $\Wc_c (\g_N)$ is topologically free.

\begin{rem}\label{ree1}
  The algebra $\Y_0  (\g_N)$ at $h=1$ becomes the (ordinary) twisted Yangian associated with the Lie algebra $\g_N$ over the complex field; cf. \cite[Sect. 3]{MNO} or  \cite[Ch. 2]{M}.
\end{rem}

Consider the degree operator on   $\Y^+ (\gl_N)$ given by $\deg t_{ij}^{(-r)}=-r$. It defines the ascending filtration over the dual Yangian such that the corresponding graded algebra $\grrr \Y^+ (\gl_N)$
is isomorphic to the universal enveloping algebra $U(\gl_N \ot t^{-1}\CC[t^{-1}])\ot_\CC \CC[[h]]$. The corresponding isomorphism is defined by the  assignments
\beq\label{cmp1}
\bar{t}_{ij}^{(-r)}\mapsto e_{ij}\ot t^{-r},
\eeq
 where $\bar{t}_{ij}^{(-r)}$ denote the images of the elements  $t_{ij}^{(-r)}$ in the corresponding component of $\grrr \Y^+ (\gl_N)$; see \cite[Sect. 2.2]{JKMY}  for more details.

Relations \eqref{srel1} imply that the algebra $\Y^+(\g_N)$  is generated by the elements
\beq\label{gg1}
s_{ij}^{(-2r)},\quad i>j\Fand
s_{ij}^{(-2r+1)},\quad i\geqslant j
\qquad\quad\text{for}\quad\qquad
r=1,2,\ldots
\eeq
in the symmetric case  $\g_N =\oo_N$, and by the elements
\beq\label{gg2}
s_{ij}^{(-2r)},\quad i\geqslant j\Fand
s_{ij}^{(-2r+1)},\quad i> j
\qquad\quad\text{for}\quad\qquad
r=1,2,\ldots
\eeq
in the antisymmetric case $\g_N=\spp_N$. 
Let us   write $\bar{s}_{ij}^{(-r)}$ for the images of the elements $s_{ij}^{(-r)}$ in the corresponding graded algebra $\grrr \Y^+ (\g_N)$. Using \eqref{Sseries} we can compute  the images of  the generators   \eqref{gg1} and \eqref{gg2}  under the isomorphism \eqref{cmp1}, thus getting
\beq\label{img64}
\bar{s}_{ij}^{(-r)}\mapsto 
\sum_{k=1}^N  g_{kj}\ts e_{ik}\ot t^{-r}
+
(-1)^{r-1} \sum_{k=1}^N g_{ik}\ts e_{jk}\ot t^{-r}
.
\eeq

Consider the twisted polynomial current Lie algebra
$$
(\gl_N \ot t^{-1}\CC[t^{-1}])^\sigma
=\left\{
A(t)\in \gl_N \ot t^{-1}\CC[t^{-1}]\,:\,
\sigma(A(t))=A(-t)
\right\},
$$
where the involutive automorphism $\sigma$ is given by \eqref{sigma}.
Its elements are  polynomials in $t^{-1}$ of the form $\sum_{k<0} A_k \ot t^k$ such that their even (resp. odd) coefficients $A_{2i}$ (resp. $A_{2i+1}$) belong to the subalgebra $\g_N$ (resp. (-1)-eigenspace of $\sigma$).
Thus,  the Lie algebra $(\gl_N \ot t^{-1}\CC[t^{-1}])^\sigma$ is spanned by the images of
$\bar{s}_{ij}^{(-r)}$  given by \eqref{img64}, so that we have

\begin{pro}\label{prro} 
The restriction of \eqref{cmp1} gives the isomorphism between the graded algebra $\grrr \Y^+ (\g_N)$ and the universal enveloping algebra $U((\gl_N \ot t^{-1}\CC[t^{-1}])^\sigma)\ot_\CC \CC[[h]]$.
\end{pro}

 The analogous results can be also obtained for the generators of the twisted  Yangian; see \cite[Sect. 3]{MNO} or \cite[Ch. 2]{M}. In fact,  Proposition \ref{prro} was established by closely following the approach from these references.

\subsection{Vacuum module \texorpdfstring{$\Wc_c (\g_N)$}{Mc(gN)} as a quasi \texorpdfstring{$\Vc_{2c} (\gl_N)$}{V2c(glN)}-module}\label{sec22}
We   now present the main result of this section, the construction of the structure of quasi $\Vc_{2c} (\gl_N)$-module  over $\Wc_c\coloneqq \Wc_c (\g_N) $. More precisely, following the   definition of   {\em quasi module} for $h$-adic nonlocal vertex algebra \cite[Def. 2.23]{Li},   we  construct the $\CC[[h]]$-module map
\begin{align*}
Y_{\Wc_c  }(z)\colon \Vc_{2c} (\gl_N)\ot \Wc_c (\g_N)&\to \Wc_c (\g_N)((z))_h\\
v\ot w&\mapsto Y_{\Wc_c  }(z)(v\ot w)=Y_{\Wc_c  }(v,z)w=\sum_{r\in\mathbb{Z}} v_r w\ts  z^{-r-1}
\end{align*}
which satisfies 
$
Y_{\Wc_c  }(\vac,z)w=w$  for all $w\in  \Wc_c (\g_N) $ 
and possesses the {\em quasi weak associativity} property:
For any positive integer $k $ and elements $u,v\in \Vc_{2c} (\gl_N)$, $w\in \Wc_c (\g_N)$
there exists a nonzero polynomial $p(x_1,x_2)\in\CC[x_1,x_2]$
such that
\begin{align}
&p(z_0 +z_2,z_2)  Y_{\Wc_c  }(u,z_0 +z_2)Y_{\Wc_c  }(v,z_2)  w \non\\
&\qquad\qquad -   p(z_0 +z_2,z_2)Y_{\Wc_c }\big(Y(u,z_0)v,z_2\big)w\,\,
\in\,\, h^k  \Wc_c (\g_N) [[z_0^{\pm 1},z_2^{\pm 1}]].\label{mod2}
\end{align}
 
In order to present the aforementioned result, we need the following notation for the operators on $(\ndo\CC^N)^{\ot n}\ot \Wc_c(\gl_N)$, where  $u=(u_1,\ldots ,u_n)$ and $v=(v_1,\ldots ,v_m)$,
\begin{align}
&S_{[n]}^{+}(u|z)=\prod_{i=1,\ldots,n}^{\longrightarrow} \left(S_{i}^+ (z+u_i)   \wvr{R}_{i\, i+1}^t(-2z-u_i-u_{i+1})\ldots \wvr{R}_{in}^t(-2z-u_i-u_n)  \right),\label{esenplus}\\
&S_{[n]}(u|z)=\prod_{i=1,\ldots,n}^{\longrightarrow} \left(S_{i} (z+u_i)    \wvr{R}_{i\, i+1}^t (-2z-u_i-u_{i+1}-hc)\ldots \wvr{R}_{in}^t (-2z-u_i-u_n-hc)   \right).\non
\end{align}
As before,  $S_{[n]}^{+}(u )$ and $S_{[n]}(u)$ stand for the analogous expressions with the variable $z$ omitted; recall \eqref{omt}.
Using   Yang--Baxter equation \eqref{ybe} and relations \eqref{usflg}--\eqref{RSRS3}
one easily verifies the following equalities    
for   operators on
$(\ndo\CC^N)^{\ot n}\ot(\ndo\CC^N)^{\ot m}\ot \Wc_c(\g_N)$:
\begin{align}
&R_{nm}^{12}(u|v|z-w)\ts S_{[n]}^{+13}(u|z)\ts \wndr{R}_{nm}^{t\,12}(-u|-v|-z-w)\ts S_{[m]}^{+23}(v|w)\non\\
&\,=S_{[m]}^{+23}(v|w)\ts \wndr{R}_{nm}^{t\,12}(-u|-v|-z-w)\ts S_{[n]}^{+13}(u|z)\ts R_{nm}^{12}(u|v|z-w),\label{rsrs12}\\
&R_{nm}^{12}(u|v|z-w)\ts S_{[n]}^{13}(u|z)\ts \wndr{R}_{nm}^{t\, 12}(-u|-v|-z-w-hc)\ts S_{[m]}^{23}(v|w)\non\\
&\,=S_{[m]}^{23}(v|w)\ts \wndr{R}_{nm}^{t\, 12}(-u|-v|-z-w-hc)\ts S_{[n]}^{13}(u|z)\ts R_{nm}^{12}(u|v|z-w),\label{rsrs22}\\
&\wvr{R}_{nm}^{12}(u|v|z-w+3hc/2)\ts S_{[n]}^{13}(u|z)\ts \wndr{\wvr{R}}_{nm}^{t\,12}(-u|-v|-z-w+hc/2)\ts S_{[m]}^{+23}(v|w)\non\\
&\,=S_{[m]}^{+23}(v|w)\ts \wndr{\wvr{R}}_{nm}^{t\,12}(-u|-v|-z-w-3hc/2)\ts S_{[n]}^{13}(u|z)\ts \wvr{R}_{nm}^{12}(u|v|z-w-hc/2),\label{rsrs32}
\end{align}
where the meaning of superscripts $1,2,3$ is the same as in \eqref{fctr}.
The next theorem is the main result in this section. It is proved in Subsection \ref{sbsprf} below.

\begin{thm}\label{qmain}
For any $c\in\CC$ there exists a unique structure of quasi $\Vc_{2c}(\gl_N)$-module  on the vacuum module $\Wc_c(\g_N)$ such that
\beq\label{Y}
Y_{\Wc_c }(T_{[n]}^+(u)\vac,z)=
S_{[n]}^{+}(u|z)\ts S_{[n]}(u|z+hc/2)^{-1} .
\eeq
\end{thm}

\begin{rem}
The expression  for the module  map in \eqref{Y} is motivated by the {\em quantum current commutation relation} from \cite{RS}. In this particular case it takes the following form: For any integer $n\geqslant 1$ there exists an integer $r\geqslant 0$ such that the equality
\begin{align*}
&(u^2-v^2)^r\ts
\wvr{R}(u-v)\ts \mathcal{L}_1(u)\ts\wvr{R}(u-v+2hc)^{-1} \mathcal{L}_2(v)\\
&\qquad =(u^2-v^2)^r\ts\mathcal{L}_2(v)\ts\wvr{R}(v-u+2hc)^{-1} \mathcal{L}_1(u)\ts\wvr{R}(v-u)
\end{align*}
holds modulo $h^n$ with $\mathcal{L}(u)=Y_{\Wc_c }(T_{[n]}^+(0)\vac,u)$.
\end{rem}

We now derive a simple consequence of the proof of Theorem \ref{qmain}. An $\Ac_c(\g_N)$-module  $M$ is said to be {\em restricted}   if it is topologically free as a $\CC[[h]]$-module and the action of $S(z)$ on $M$ belongs to $ \ndo\CC^N \ot \om (M,M[z^{-1}]_h)$. The argument from the proof of Lemma \ref{lemma24} shows that $\Wc_c(\g_N)$ is a restricted module. One can  easily prove  that on any restricted $\Ac_c(\g_N)$-module $M$  for any $m\geqslant 1$ we have
$$
S_{[m]}(u|z)  \in (\ndo\CC^N)^{\ot m
} \ot \om (M,M[z^{-1}][[u_1,\ldots ,u_m]]_h).
$$
This is due to the fact that for any $n\geqslant 1$ the coefficients of the powers of $u$ in $\wvr{R}^t(-2z-u)$ and $S(z+u) w$, where $w\in M$,   are polynomials in $z^{-1}$ when regarded modulo $h^n$.
Thus the proof of the next corollary goes in parallel with the proof of Theorem \ref{qmain}.

\begin{kor}
Let $M$ be a restricted $\Ac_c(\g_N)$-module. There exists a unique structure of quasi $\Vc_{2c}(\gl_N)$-module  on   $M$ such that
$$
Y_{M}(T_{[n]}^+(u)\vac,z)=
S_{[n]}^{+}(u|z)\ts S_{[n]}(u|z+hc/2)^{-1} .
$$
\end{kor}

\subsection{Proof of Theorem \ref{qmain}}\label{sbsprf}
The proof of  Theorem \ref{qmain}  is divided in three  lemmas which verify the requirements imposed by the definition of quasi module, as given in Subsection \ref{sec22}. 
Even though their proofs  go in parallel with the proof of  \cite[Thm. 2.7]{K19}, we provide most of the details in order to take care of the variations  which occur  in this setting. 
Throughout the proof, we  use the ordered product notation, where the subscript of the product symbol determines the order of the tensor factors. More specifically, for any   elements $a=\sum_i a_1^{(i)}\ot a_2^{(i)}$ and $b=\sum_j b_1^{(j)} \ot b_2^{(j)}$  of $\ndo \CC^N \ot \ndo \CC^N$ we write
$$
a\cdotlr b =\sum_{i,j} a_1^{(i)}b_1^{(j)}\ot  b_2^{(j)} a_2^{(i)}\Fand
a\cdotrl b =\sum_{i,j} b_1^{(j)}a_1^{(i)}\ot   a_2^{(i)}b_2^{(j)}.
$$
Using this notation we can express the crossing symmetry property \eqref{csym} as
\beq\label{csym2}
\wvr{R}(u)^{-1}\cdotlr \wvr{R}(u+hN)=\wvr{R}(u)^{-1}\cdotrl \wvr{R}(u+hN)=1.
\eeq

\begin{lem}
Formula \eqref{Y}, together with $Y_{\Wc_c }( \vac,z)=1_{\Wc_c(\g_N)}$, defines a unique $\CC[[h]]$-module map
$\Vc_{2c} (\gl_N)\ot \Wc_c (\g_N) \to \Wc_c (\g_N)[[z^{\pm 1}]]$.
\end{lem}

\begin{prf}
The coefficients of the matrix entries of  $T_{[n]}^+(u )\vac$ with $n\geqslant 0$ span an $h$-adically dense $\mathbb{C}[[h]]$-submodule of $\Vc_{2c}(\gl_N)$, so \eqref{Y} uniquely determines the quasi module map. Hence it remains to verify that the   quasi module map is well-defined by \eqref{Y}. It is sufficient to show that it maps
 the ideal of  relations \eqref{RTT1}   to itself. Indeed, by the Poincar\'e--Birkhoff--Witt theorem for the double Yangian   \cite[Thm. 2.2]{JKMY}, the dual Yangian coincides with the algebra defined by the generators $t_{ij}^{(-r)}$, where $r\geqslant 1$ and $i,j=1,\ldots ,N$, subject to the defining relations given by \eqref{RTT1}. 

First, we introduce some notation and establish some Yang--Baxter-like identities which we shall use in the later stage of the proof.
Let $n\geqslant 2$ be an arbitrary integer    and $u=(u_1,\ldots ,u_n)$ a family of variables.  For all   $i,j=1,\ldots ,n$ such that $i\neq j$ we write
$$
R_{ij}=R_{ij}(u_i -u_j),\quad \wvr{R}_{ij}=\wvr{R}^t_{ij}(-2z-u_i-u_j),\quad
\wht{R}_{ij} =\wvr{R}_{ij}^t(-2z-u_i-u_j-2hc) .
$$
 The Yang--Baxter equation \eqref{ybe} and the unitarity property \eqref{uni} imply the following identities for all  $1\leqslant j<k<k+1<l\leqslant n$:
\begin{gather}
R_{k\tss k+1} \wvr{R}_{jk} \wvr{R}_{j\tss k+1} =\wvr{R}_{j\tss k+1}   \wvr{R}_{jk}  R_{k\tss k+1},\quad
R_{k\tss k+1} \wvr{R}_{kl} \wvr{R}_{k+1\tss l} =\wvr{R}_{k+1\tss l}   \wvr{R}_{kl}  R_{k\tss k+1},\label{prf2}\\
   R_{k\tss k+1} \wht{R}_{k+1\tss l}  ^{-1}  \wht{R}_{kl}  ^{-1} = \wht{R}_{kl}   ^{-1} \wht{R}_{k+1\tss l}  ^{-1} R_{k\tss k+1},\quad
   R_{k\tss k+1} \wht{R}_{j\tss k+1}  ^{-1} \wht{R}_{jk}  ^{-1} = \wht{R}_{jk}  ^{-1} \wht{R}_{j\tss k+1}  ^{-1} R_{k\tss k+1}.\label{prfx}
\end{gather}

For any $k=1,\ldots , n-1$
defining relation \eqref{RTT1} for the dual Yangian implies 
\beq\label{prf0}
R_{k\tss k+1}T_{[n]}^+(u )=T_{1}^+(u_{1})\ldots T_{k-1}^+(u_{k-1})T_{k+1}^+(u_{k+1})T_{k}^+(u_{k})T_{k+2}^+(u_{k+2})\ldots T_{n}^+(u_{n})R_{k\tss k+1}.
\eeq
We shall prove that the   difference of the left and the right-hand side of \eqref{prf0} belongs to the kernel of $Y_{\Wc_c }(\cdot, z)$, which implies the lemma.    It is clear from \eqref{Y} that the image of the left-hand side under  $Y_{\Wc_c }(\cdot, z)$ is equal to
\beq\label{exprr1}
R_{k\tss k+1}\ts S_{[n]}^{+}(u|z)\ts S_{[n]}(u|z+hc/2)^{-1} .
\eeq
As for the right-hand side, using \eqref{Y} again, we find that its image is given by
\beq\label{exprr2}
P_{k\ts k+1}\ts S_{[n]}^{+}(u^{(k)}|z)\ts S_{[n]}(u^{(k)}|z+hc/2)^{-1}\ts P_{k\ts k+1} \ts R_{k\tss k+1},
\eeq
where $u^{(k)}=(u_1,\ldots ,u_{k-1},u_{k+1},u_k,u_{k+2},\ldots ,u_n)$ and $P_{k\ts k+1}$ stands for the action of the permutation operator $P$    on the tensor factors $k$ and $k+1$; recall \eqref{perm}. 
By employing the relation \eqref{RSRS1} and the equalities in \eqref{prf2} we find
\beq\label{exprr3}
R_{k\tss k+1}\ts S_{[n]}^{+}(u|z)=P_{k\ts k+1}\ts S_{[n]}^{+}(u^{(k)}|z)\ts P_{k\ts k+1}\ts R_{k\tss k+1}.
\eeq
Analogously, by using \eqref{RSRS2} and the equalities in \eqref{prfx} we obtain
\beq\label{exprr4}
R_{k\tss k+1}\ts S_{[n]}(u|z+hc/2)^{-1}=P_{k\ts k+1}\ts S_{[n]}(u^{(k)}|z+hc/2)^{-1}\ts P_{k\ts k+1}\ts R_{k\tss k+1}.
\eeq
Finally,  from  \eqref{exprr3} and \eqref{exprr4} we easily see that the expressions in \eqref{exprr1} and \eqref{exprr2} coincide, so that the map $Y_{\Wc_c }(\cdot, z)$ is well-defined by \eqref{Y}, as required.
\end{prf}

\begin{lem}\label{lemma24}
For any $v\in \Vc_{2c} (\gl_N)$ and $w\in \Wc_c (\g_N)$ the series $Y_{\Wc_c }(v,z)w$ belongs to $\Wc_c (\g_N)((z))_h$, i.e. it possesses   finitely many negative powers of $z$ modulo $h^k$ for all $k\geqslant 0$.
\end{lem}

\begin{prf}
Note that the coefficients of the matrix entries of  $S_1^+(v_1)\ldots S_m^+(v_m)$ with $m\geqslant 0$ span an $h$-adically dense $\mathbb{C}[[h]]$-submodule of $\Wc_c (\g_N)$. Furthermore, from
\beq\label{n76}
S_{[m]}^{+}(v )=\prod_{i=1,\ldots,m}^{\longrightarrow} \left(S_{i}^+ ( v_i)\ts  \wvr{R}_{i\, i+1}^t( -v_i-v_{i+1})\ldots \wvr{R}_{im}^t( -v_i-v_m)  \right)
\eeq
one easily shows that the coefficients of the matrix entries of  $S_{[m]}^{+}(v )$ with $m\geqslant 0$ span an $h$-adically dense $\mathbb{C}[[h]]$-submodule of $\Wc_c (\g_N)$. Indeed, we can move all $R$-matrices to the left-hand side of \eqref{n76} by using the identity
$
\wvr{R}^t( u)\cdotrl\wvr{R}^t( -u)=1
$; recall \eqref{uni}.
Thus,  on the right-hand side we get $S_1^+(v_1)\ldots S_m^+(v_m)$, while   the left-hand side is expressed  in terms of $S_{[m]}^{+}(v )$. Hence the coefficients of the matrix entries of  $S_1^+(v_1)\ldots S_m^+(v_m)$ can be expressed in terms of the coefficients of the matrix entries of $S_{[m]}^{+}(v )$, as required.

By the preceding discussion, it is sufficient to check that for any integers $m,n\geqslant 1$ the image of \eqref{n76}  under  \eqref{Y} belongs to $\Wc_c (\g_N)((z))_h$. Clearly, this image equals
\beq\label{nm4}
S_{[n]}^{+13}(u|z)\ts  S_{[n]}^{13}(u|z+hc/2)^{-1}S_{[m]}^{+23}(v ).
\eeq
Using relation \eqref{rsrs32}, along with the crossing symmetry property   \eqref{csym2}, we rewrite it as
\begin{align}
&S_{[n]}^{+13}(u|z)\bigg(\wvr{R}_{nm}^{12}(u|v|z +2hc +hN)\cdotrl\bigg( \wndr{\wvr{R}}_{nm}^{t\,12}(-u|-v|-z  )\ts  S_{[m]}^{+23}(v ) \bigg.\bigg.\non\\
\times &\wvr{R}_{nm}^{12}(u|v|z  )^{-1}S_{[n]}^{13}(u|z+hc/2)^{-1}\vac \wndr{\wvr{R}}_{nm}^{t\,12}(-u|-v|-z -2hc )^{-1} \bigg.\bigg.\bigg)\bigg) .\label{nm5}
\end{align}
Note that by \eqref{Sseries} we have $S(x)^{\pm 1}\vac=G^{\pm}\ot \vac$. In addition, recall that  the normalized $R$-matrix   $\wvr{R}(x) $ belongs to $(\ndo\CC^N)^{\ot 2}  [x^{-1}]_h$. Therefore, we   conclude that for any choice of    $a_1,\ldots ,a_n,b_1,\ldots ,b_m,k$, by regarding the expression in \eqref{nm5} modulo
\beq\label{monomii}
u_1^{a_1}\ldots u_n^{a_n}v_1^{b_1}\ldots v_m^{b_m} h^k,
\eeq
 we obtain only finitely many negative powers of the variable $z$. Indeed, the negative powers of $z$ in \eqref{nm5}   come  from the $R$-matrices only. However,   the expression contains finitely many $R$-matrices and   each of them produces   finitely many negative powers of  $z$ modulo \eqref{monomii}.
Hence  the image of $Y_{\Wc_c}(\cdot, z)$ belongs to $\Wc_c (\g_N)((z))_h $, as required.
\end{prf}

\begin{lem}
The map $Y_{\Wc_c }(\cdot,z)$ possesses the quasi weak associativity property \eqref{mod2}.
\end{lem}

\begin{prf}
Consider  the second summand in \eqref{mod2}.
Let us apply  the vertex operator map $Y(\cdot,z_0)$ of $\Vc_{2c}(\gl_N)$, as defined in  \eqref{qva1}, on the series
\beq\label{as1}
T_{[n]}^{+13}(u)\ts   T_{[m]}^{+24}(v)\ts  (\vac\ot \vac)\in(\ndo\CC^N)^{\ot n}\ot (\ndo\CC^N)^{\ot m}\ot \Vc_{2c}(\gl_N)\ot \Vc_{2c}(\gl_N)[[u,v]],
\eeq 
where $u=(u_1,\ldots ,u_n)$ and $v=(v_1,\ldots ,v_m)$.
By using the relation \eqref{rtt3} at the level $2c$, the crossing symmetry   \eqref{csym2} and the identity $T_{[n]}^{13}(u|z_0 +hc)^{-1}\vac=1^{\ot n}\ot\vac$, we obtain
\begin{align}
&\,T_{[n]}^{+13}(u|z_0)\ts  T_{[n]}^{ 13}(u|z_0+hc)^{-1}  T_{[m]}^{+23}(v)\vac
=\wvr{R}_{nm}^{12,+} \cdotrl \left(T_{[n]}^{+13}(u|z_0)\ts   T_{[m]}^{+23}(v)\ts  \wvr{R}_{nm}^{12,-}   \right)
,\label{ztr1}
\end{align}
where the terms $\wvr{R}_{nm}^{12,\pm}$ are defined by
$$
\wvr{R}_{nm}^{12,+}=\wvr{R}_{nm}^{12}(u|v|z_0+2hc+hN)
\fand
\wvr{R}_{nm}^{12,-}=\wvr{R}_{nm}^{12}(u|v|z_0)^{-1} .
$$
Next, we apply the map $Y_{\Wc_c }(\cdot,z_2)$ to the right-hand side of \eqref{ztr1}. By \eqref{Y} we get
\begin{align}
&\wvr{R}_{nm}^{12,+} \cdotrl \left(   S_{[n]}^{+13}(u|z_2+z_0)\ts  \wvr{\wndr{R}}_{nm}^{12,+}\ts   S_{[m]}^{+23}(v|z_2)\right.\non\\
& \hspace{40pt}\left.\times S_{[m]}^{ 23}(v|z_2+hc/2)^{-1}\ts  
\wvr{\wndr{R}}_{nm}^{12,-}\ts  
S_{[n]}^{ 13}(u|z_2+z_0+hc/2)^{-1}\ts  \wvr{R}_{nm}^{12,-}  \right),\label{tz6}
\end{align}
where the terms $\wvr{\wndr{R}}_{nm}^{12,\pm}$ are given by
$$
\wvr{\wndr{R}}_{nm}^{12,+}
=\wvr{\wndr{R}}_{nm}^{t\, 12}(-u|-v|-2z_2-z_0)
\fand
\wvr{\wndr{R}}_{nm}^{12,-}
=\wvr{\wndr{R}}_{nm}^{t\, 12}(-u|-v|-2z_2-z_0-2hc)^{-1}.
$$
Using    relation \eqref{rsrs22}   we rearrange the last four factors in \eqref{tz6}, thus getting
\begin{align*}
&\wvr{R}_{nm}^{12,+} \cdotrl \left(S_{[n]}^{+13}(u|z_2+z_0)\ts  \wvr{\wndr{R}}_{nm}^{12,+} \ts  S_{[m]}^{+23}(v|z_2)\right.\non\\
&\hspace{40pt}\left.\times \wvr{R}_{nm}^{12,-}\ts   
S_{[n]}^{ 13}(u|z_2+z_0+hc/2)^{-1} \ts  
\wvr{\wndr{R}}_{nm}^{12,-}\ts  
S_{[m]}^{ 23}(v|z_2+hc/2)^{-1}\right)
. 
\end{align*}
Finally, we employ relation \eqref{rsrs32} to rewrite this as
\begin{align*}
&\wvr{R}_{nm}^{12,+}    \cdotrl \left(S_{[n]}^{+13}(u|z_2+z_0)\ts  
S_{[n]}^{ 13}(u|z_2+z_0+hc/2)^{-1} 
\right.\\
&\hspace{40pt}\left.\times \wvr{R}_{nm}^{12 } (u|v|z_0+2hc)^{-1}\ts   S_{[m]}^{+23}(v|z_2)\ts  
 (\wvr{\wndr{R}}_{nm}^{12,-})^{-1}\ts  
\wvr{\wndr{R}}_{nm}^{12,-}\ts  
S_{[m]}^{ 23}(v|z_2+hc/2)^{-1}\right)
. 
\end{align*}
It remains to observe that by the crossing symmetry property \eqref{csym2}  the $R$-matrix factors in the given expression cancel, so that it equals
 \begin{align}
 S_{[n]}^{+13}(u|z_2+z_0)\ts  
S_{[n]}^{ 13}(u|z_2+z_0+hc/2)^{-1} 
   S_{[m]}^{+23}(v|z_2)\ts  
S_{[m]}^{ 23}(v|z_2+hc/2)^{-1} \label{tz7}
. 
\end{align}

Recall that the expression in \eqref{tz7} corresponds to the second summand in \eqref{mod2}. As for the first summand, one immediately sees from \eqref{Y} that it is equal to
 \begin{align}
 S_{[n]}^{+13}(u|z_0+z_2)\ts  
S_{[n]}^{ 13}(u|z_0+z_2+hc/2)^{-1} 
   S_{[m]}^{+23}(v|z_2)\ts  
S_{[m]}^{ 23}(v|z_2+hc/2)^{-1} \label{tz8}
. 
\end{align}
Observe that  \eqref{tz7} and \eqref{tz8} do not coincide as the former should be expanded in negative powers of $z_2$ and the latter in the negative powers of $z_0$. However, let us apply   both \eqref{tz7} and \eqref{tz8}  on an arbitrary element $w\in  \Wc_c (\g_N)$. 
By arguing as  in the proof of Lemma \ref{lemma24}, one   checks
that there exists $r\geqslant 0$ such that the  products of the 
resulting expressions with $p(z_0+z_2,z_2)$, where
$p(x_1,y_1)=x_1^r (x_1+x_2)^r$,
coincide modulo \eqref{monomii}. Thus the map $Y_{\Wc_c }(\cdot,z)$ satisfies the quasi weak associativity  requirement  \eqref{mod2}.\footnote{Note that    the polynomial $p(x_1,y_1)=x_1^r (x_1+x_2)^r$ can not be replaced by the simpler polynomial $q(x_1,y_1)=x_1^r$. In other words, the quasi module map 
$Y_{\Wc_c }(\cdot,z)$ does not need to satisfy the usual {\em weak associativity}, which is obtained from \eqref{mod2} by replacing
$p(z_0 +z_2,z_2)$ with $(z_0+z_2)^r$.  }
\end{prf}

\section{Central elements in  \texorpdfstring{$\wtld{\Ac}_{c}(\g_N)$}{Ac(gN)} and invariants in \texorpdfstring{$\Wc_c(\g_N)$}{Mc(gN)}}\label{sec3}
In this section, we construct families of central elements of the  completed algebras $\wtld{\Ac}_{c}(\g_N)$   and invariants of the quasi module $\Wc_c(\g_N)$. Our constructions employ  the quasi module map established by Theorem \ref{qmain}. Moreover, they rely on the fusion procedure for the Yang $R$-matrix and the explicit formulae for certain families of central elements of the quantum    vertex algebra $\Vc_{c}(\gl_N)$, so we start by briefly recalling these results.

\subsection{Fusion procedure and the  center of \texorpdfstring{$\Vc_{c}(\gl_N)$}{Vc(glN)}}
Let us recall the fusion procedure for the Yang $R$-matrix \eqref{yang} found in \cite{J}; see also \cite[Sect. 6.4]{M} for more details. Let $\nu$ be a Young diagram with $n$ boxes of length   less than or equal to $N$   and let $\Uc$ be a standard $\nu$-tableau with entries in $\left\{1,\ldots,n\right\}$. For any $p=1,\ldots,n$   the contents $c_p=c_p(\Uc)$ of $\Uc$ are defined by   $c_p = j-i$ if $p$ occupies the box $(i,j)$ of $\Uc$.
Let $e_{\Uc}$ be  the primitive idempotent in the group algebra $\CC[\Sym_n]$ of the symmetric group $\Sym_n$ associated with $\Uc$
through the use of the orthonormal Young bases for the irreducible representations of $\Sym_n$.
The  group $\Sym_n$ acts on   $(\CC^N)^{\ot\ts n}$ by permuting the tensor
factors. Let $\Ec^{}_{\Uc}$ be the image
of $e^{}_{\Uc}$ with respect to this action.
By \cite{J}, the consecutive evaluations $u_k =hc_k$  
of 
$$R(u_1,\dots,u_n)=\prod^{\longrightarrow }_{1\leqslant i< j\leqslant n} R_{ij}(u_i -u_j) $$
 are well-defined. Furthermore,
  the result is proportional to $\Ec^{}_{\Uc}$, i.e. we have
\beq\label{fusion}
R(u_1,\dots,u_n)\big|_{u_1=hc_1}
\big|_{u_2=hc_2}\dots \big|_{u_n=hc_n}=\pi_\nu\ts \Ec^{}_{\Uc},
\eeq
where $\pi_\nu$ stands for the product of all hook lengths of the boxes of $\nu$.

Consider the quantum affine vertex algebra at the critical level  $\Vc_{-N}(\gl_N)$. 
Let
\beq\label{variable}
u_\nu =(u_1,\ldots,u_n),\qquad \text{where} \qquad
u_k = u+hc_k \quad\text{for}\quad k=1,\ldots,n.
\eeq 
Due to \cite[Thm. 2.4]{JKMY}, all  coefficients of   the series 
$$
\TT_{\nu}^+ (u )=\tr_{1,\ldots,n}\ts \Ec_{\Uc}\ts  T_{[n]}^+ (u_\nu) \vac=\tr_{1,\ldots,n}\ts \Ec_{\Uc}\ts  T_1^+ (u_1)\ldots T_n^{+}(u_n) \vac\in \Vc_{-N}(\gl_N)[[u]],
$$
 where the trace   is taken over all $n$ copies of $\ndo\CC^N$,
belong to the center $\z(\Vc_{-N}(\gl_N))$ of the quantum vertex algebra $\Vc_{-N}(\gl_N)$.
As for the noncritical level  $c\neq -N$, all coefficients of the {\em quantum determinant}
\begin{align}\label{kvdet}
 \qdet T^+(u)=\sum_{\sigma\in\Sym_N}\sgn \sigma\cdot
t^+_{ \sigma(1)\tss 1}(u) \dots t^+_{ \sigma(N)\tss N}(u-(N-1)h)\vac\in \Vc_{c}(\gl_N)[[u]] 
\end{align}
belong to the center $\z(\Vc_{c}(\gl_N))$ of the quantum vertex algebra $\Vc_{c}(\gl_N)$; see \cite[Prop. 4.10]{JKMY}.

\subsection{Central elements of the completed algebra  \texorpdfstring{$\wtld{\Ac}_{c}(\g_N)$}{Ac(gN)}}

Let $I_p$ for $p\geqslant 1$ be the left ideal of the double Yangian $\DY_c(\gl_N)$ at the level 
 $c\in\CC$  generated by all elements $t_{ij}^{(r)}$ such that $r\geqslant p$.  
 Following \cite{JKMY}, we define the completed  double Yangian $\wtld{\DY}_c(\gl_N)$ at the level 
  $c$    as  the $h$-adic  completion of the inverse limit
 $\lim_{\longleftarrow} \,  \DY_c(\gl_N)/  I_p.
 $
Finally, we introduce the algebra    $\wtld{\Ac}_c(\g_N)$  as  the $h$-adic  completion of the inverse limit
$$\lim_{\longleftarrow} \,  \Ac_c(\g_N)/ (\Ac_c(\g_N)\cap I_p).
$$

\subsubsection{Critical level}
In this subsection, we consider the algebra $\wtld{\Ac}_{-N/2}(\g_N)$. For any integer $n=1,\ldots ,N$ introduce  the Laurent series with coefficients in $\wtld{\Ac}_{-N/2}(\g_N)$ by
\begin{align}
& \sss_{\nu } (u)=\tr_{1,\ldots,n}\ts   \Ec_\Uc \ts S^+_{[n]}(u_\nu)\ts S_{[n]}(u_\nu-hN/4)^{-1},\label{6c1}
\end{align}
where the family of variables $u_\nu$ is given by \eqref{variable}. It is worth noting that the action of the series $\sss_{\nu } (u)$ on the vacuum module $\Wc_{-N/2}(\g_N)$ is given by
\beq\label{sfh1}
\sss_{\nu } (u)=Y_{\Wc_{-N/2}}(\TT_{\nu}^+(0),u).
\eeq

In what follows, we shall often use the arrow at the top of the symbol to indicate  that the   factors are arranged in the opposite order, e.g., 
for $\wvr{R}_{i\ts n+1}=\wvr{R}_{i\ts n+1}(u_i-v)$ we have
$$
  \wvr{R}_{n1} (u |v)=  \wvr{R}_{1\ts n+1} \ldots \wvr{R}_{n\ts n+1} 
	\Fand
  \cev{\wvr{R}}_{n1} (u |v)=\wvr{R}_{n\ts n+1} \ldots \wvr{R}_{1\ts n+1} .
$$
In the next lemma $u_0$ denotes a single variable while $u_\nu$ again stands for the family \eqref{variable}.
\begin{lem}\label{lk31}
The following equalities hold on $\ndo\CC^N \ot (\ndo\CC^N)^{\ot n}\ot \wtld{\Ac}_{-N/2}(\g_N)$:
\begin{align}
&\Ec_{\Uc} \wvr{R}_{1n}^{01}(u_0|u_\nu)= \cev{\wvr{R}}_{1n}^{01}(u_0|u_\nu)\Ec_{\Uc},\qquad \Ec_{\Uc} \wvr{R}_{1n}^{01}(u_0|u_\nu)^{-1}= \cev{\wvr{R}}_{1n}^{01}(u_0|u_\nu)^{-1}\Ec_{\Uc},\label{xxx1}\\
&\Ec_{\Uc} \cev{\wvr{\wndr{R}}}_{1n}^{t\,01}(u_0|u_\nu)=  \wvr{\wndr{R}} _{1n}^{t\,01}(u_0|u_\nu)\Ec_{\Uc},\,\,\,\,  \Ec_{\Uc} \wvr{\wndr{R}}_{1n}^{t\,01}(-u_0|-u_\nu) = \cev{\wvr{\wndr{R}}}_{1n}^{t\,01}(-u_0|-u_\nu) \Ec_{\Uc},\label{xxx3}\\ 
&\Ec_{\Uc}S^{+12}_{[n]}(u_\nu)=\cev{S}^{+12}_{[n]}(u_\nu)\Ec_{\Uc},\qquad \Ec_{\Uc}S_{[n]}^{12}(u_\nu-hN/4)^{-1}=\cev{S}^{12}_{[n]}(u_\nu-hN/4)^{-1}\Ec_{\Uc},\label{xxx5}
\end{align}
where the $n+1$ copies
of $\ndo\CC^N$ are  labeled by $0,\ldots,n$,
 $\Ec_{\Uc}$  is applied on the tensor factors $1,\ldots,n$ and the superscripts indicate the tensor factors as follows:
\beq\label{ten1ten1}
\overbrace{\ndo\CC^N}^{0} \ot \overbrace{(\ndo\CC^N)^{\ot n}}^{1}\ot \overbrace{\wtld{\Ac}_{-N/2}(\g_N)}^{2}.
\eeq
\end{lem}

\begin{prf}
The identities in \eqref{xxx1}   are proved in \cite[Lemma 3.1]{K19}.
The remaining equalities can be verified analogously, by using the Yang--Baxter equation \eqref{ybe},  the fusion procedure \eqref{fusion}  and the  relations   \eqref{rsrs12} and \eqref{rsrs22}.
\end{prf}

The next theorem  is the main result of this section. As with the similar constructions of the families of central elements in the completions of the   double Yangian and the reflection algebra, such as \cite[Thm. 4.4]{JKMY}, \cite[Thm. 4.4]{JY} and \cite[Thm. 3.2]{K19}, its proof relies on the usual $R$-matrix techniques and the fusion procedure for the Yang $R$-matrix.

\begin{thm}\label{main}
All coefficients of $\BB_{\nu} (u)$   belong to the center of the   algebra $\wtld{\Ac}_{-N/2}(\g_N)$.
\end{thm}

\begin{prf}
Let us prove 
\beq\label{c2}
 S (u_0)\ts  \BB_{\nu}  (u)=\BB_{\nu}  (u)\ts S (u_0).
\eeq
By applying $S (u_0)$ to the right-hand side of \eqref{6c1} we get
\beq\label{c3}
\tr_{1,\ldots,n}\ts  \Ec_\Uc\ts S_0(u_0)\ts S^{+12}_{[n]}(u_\nu)\ts S_{[n]}^{12}(u_\nu -hN/4)^{-1}.
\eeq
Observe that our  notation coincides with \eqref{ten1ten1} and, in particular, that
the trace is taken over the tensor factors $1,\ldots ,n$ while $S(u_0)$ is applied on the $0$-th factor.
By employing  the unitarity \eqref{uni} and the relation \eqref{rsrs32}   we express \eqref{c3} as
\begin{align}
\tr_{1,\ldots,n} &\ts \Ec_\Uc\,\bigg(
 \cev{\wndr{\wvr{R}}}_{1n}^{t\, 01  }(u_0+hN/4|u_\nu )  
 \cdotrl
\Big( \wvr{R}_{1n}^{01}(u_0-3hN/4|u_\nu)^{-1} 
S_{[n]}^{+12}(u_\nu)\Big.\bigg.\non\\
&\bigg.\Big. 
\wndr{\wvr{R}}_{1n}^{t \, 01}(-u_0+3hN/4|-u_\nu)\ts 
S_{0}(u_0)\ts 
\wvr{R}_{1n}^{01}(u_0+hN/4|u_\nu) \ts 
S_{[n]}^{12}(u_\nu-hN/4)^{-1} 
\Big)\bigg)
.\label{c5}
\end{align}
Since $ \Ec_\Uc$ is idempotent,  the first equality in \eqref{xxx3}  implies
\beq\label{c4}
\Ec_\Uc \cev{K}=  \Ec_\Uc^2 \cev{K}=\Ec_\Uc  K \Ec_\Uc=\Ec_\Uc  K \Ec_\Uc^2 \quad\text{for}\quad \cev{K}=\cev{\wndr{\wvr{R}}}_{1n}^{t\, 01  }(u_0+hN/4|u_\nu ) .
\eeq
By using \eqref{c4} we can write \eqref{c5} as
\begin{align*}
\tr_{1,\ldots,n} &\ts \Ec_\Uc\,\bigg(
 \wndr{\wvr{R}}_{1n}^{t\, 01  }(u_0+hN/4|u_\nu )  
 \cdotrl
\Big( \Ec_\Uc^2 \ts \wvr{R}_{1n}^{01}(u_0-3hN/4|u_\nu)^{-1} 
S_{[n]}^{+12}(u_\nu)\Big.\bigg.\non\\
&\bigg.\Big. 
\wndr{\wvr{R}}_{1n}^{t \, 01}(-u_0+3hN/4|-u_\nu)\ts 
S_{0}(u_0)\ts 
\wvr{R}_{1n}^{01}(u_0+hN/4|u_\nu) \ts 
S_{[n]}^{12}(u_\nu-hN/4)^{-1} 
\Big)\bigg)
. 
\end{align*}
Due to the cyclic property of the trace, this equals to
\begin{align}
\tr_{1,\ldots,n} & \ts 
   \Ec_\Uc \ts 
	\wvr{R}_{1n}^{01}(u_0-3hN/4|u_\nu)^{-1}  
S_{[n]}^{+12}(u_\nu) \ts 
\wndr{\wvr{R}}_{1n}^{t \, 01}(-u_0+3hN/4|-u_\nu)\ts 
S_{0}(u_0) \non\\
& 
\wvr{R}_{1n}^{01}(u_0+hN/4|u_\nu) \ts 
S_{[n]}^{12}(u_\nu-hN/4)^{-1}  
\Ec_\Uc\ts 
 \wndr{\wvr{R}}_{1n}^{t\, 01  }(u_0+hN/4|u_\nu ) \ts 
\Ec_\Uc 
. \label{c8}
\end{align}
By the second equality in \eqref{xxx1} and $\Ec_\Uc^2=\Ec_\Uc$, we have
$$\Ec_\Uc L=\Ec_\Uc^2 L=\Ec_\Uc \cev{L}\Ec_\Uc =\Ec_\Uc \cev{L}\Ec_\Uc^2 =\Ec_\Uc^2 L\Ec_\Uc =\Ec_\Uc L\Ec_\Uc \quad\text{for}\quad
L=\wvr{R}_{1n}^{01}(u_0-3hN/4|u_\nu)^{-1}.$$
Therefore, using the cyclic property of the trace,   we can  write \eqref{c8} as
\begin{align*}
\tr_{1,\ldots,n} & \ts 
    \wvr{R}_{1n}^{01}(u_0-3hN/4|u_\nu)^{-1}\ts 
		\Ec_\Uc\ts 
S_{[n]}^{+12}(u_\nu) \ts 
\wndr{\wvr{R}}_{1n}^{t \, 01}(-u_0+3hN/4|-u_\nu)\ts 
S_{0}(u_0) \non\\
& 
\wvr{R}_{1n}^{01}(u_0+hN/4|u_\nu) \ts 
S_{[n]}^{12}(u_\nu-hN/4)^{-1}  
\Ec_\Uc\ts 
 \wndr{\wvr{R}}_{1n}^{t\, 01  }(u_0+hN/4|u_\nu ) \ts 
\Ec_\Uc 
.  
\end{align*}
We now employ the first equality in \eqref{xxx1}, the second equality in   \eqref{xxx3} and both equalities in \eqref{xxx5}  to move the leftmost copy of $\Ec_\Uc $ to the right, which gives us 
\begin{align}
\tr_{1,\ldots,n} & \ts 
    \wvr{R}_{1n}^{01}(u_0-3hN/4|u_\nu)^{-1}  
\cev{S}_{[n]}^{+12}(u_\nu) \ts 
\cev{\wndr{\wvr{R}}}_{1n}^{t \, 01}(-u_0+3hN/4|-u_\nu)\ts 
S_{0}(u_0) \non\\
& 
\cev{\wvr{R}}_{1n}^{01}(u_0+hN/4|u_\nu) \ts 
\cev{S}_{[n]}^{12}(u_\nu-hN/4)^{-1}  
\Ec_\Uc^2 \ts 
 \wndr{\wvr{R}}_{1n}^{t\, 01  }(u_0+hN/4|u_\nu ) \ts 
\Ec_\Uc 
.  \label{c6}
\end{align}
Using \eqref{c4} and $\Ec_\Uc^2=\Ec_\Uc$ we   replace $\Ec_\Uc^2 
K \Ec_\Uc $  by $\Ec_\Uc 
\cev{K} $ in \eqref{c6}.  Next,  we use the corresponding equalities from Lemma \ref{lk31} to move the remaining copy of $\Ec_\Uc$ to the left:
\begin{align*}
\tr_{1,\ldots,n} & \ts
    \wvr{R}_{1n}^{01}(u_0-3hN/4|u_\nu)^{-1}  
		\Ec_\Uc\ts
S_{[n]}^{+12}(u_\nu) \ts
\wndr{\wvr{R}}_{1n}^{t \, 01}(-u_0+3hN/4|-u_\nu)\ts
S_{0}(u_0) \non\\
& 
\wvr{R}_{1n}^{01}(u_0+hN/4|u_\nu) 
S_{[n]}^{12}(u_\nu-hN/4)^{-1}  \ts
 \cev{\wndr{\wvr{R}}}_{1n}^{t\, 01  }(u_0+hN/4|u_\nu ) 
.   
\end{align*}
By   the relation \eqref{rsrs22}  this coincides with
$$
\tr_{1,\ldots,n}  \ts
    \wvr{R}_{1n}^{01}(u_0-3hN/4|u_\nu)^{-1} \ts
		\Ec_\Uc\ts
S_{[n]}^{+12}(u_\nu) \ts
S_{[n]}^{12}(u_\nu-hN/4)^{-1}  \ts
\wvr{R}_{1n}^{01}(u_0+hN/4|u_\nu)\ts
S_{0}(u_0).  
$$
Indeed, the   last two terms which appear on the right,  
$\wndr{\wvr{R}}_{1n}^{t \, 01}(-u_0+3hN/4|-u_\nu)$
and
$\cev{\wndr{\wvr{R}}}_{1n}^{t\, 01  }(u_0+hN/4|u_\nu ) $
cancel because of the crossing symmetry property  \eqref{csym}.
Finally, we  use    the cyclic property of the trace and \eqref{csym2}  to rewrite this as
\begin{align*}
&\tr_{1,\ldots,n}     \ts 
		\Ec_\Uc\ts 
S_{[n]}^{+12}(u_\nu) \ts 
S_{[n]}^{12}(u_\nu-hN/4)^{-1}  
\left(
\wvr{R}_{1n}^{01}(u_0-3hN/4|u_\nu)^{-1}\cdotlr
\wvr{R}_{1n}^{01}(u_0+hN/4|u_\nu)\right)
S_{0}(u_0) \\
&=\tr_{1,\ldots,n} \ts      
		\Ec_\Uc\ts 
S_{[n]}^{+12}(u_\nu) \ts 
S_{[n]}^{12}(u_\nu-hN/4)^{-1}   
S_{0}(u_0)=\BB_{\nu}  (u)\ts S (u_0),
\end{align*}
so that the commutation relation \eqref{c2}   follows.
A similar calculation verifies the equality
$$
 S^+(u_0)\ts  \BB_{\nu }  (u)= \BB_{\nu }  (u)\ts S^{+}(u_0),
$$
which,  together with \eqref{c2}, implies the statement of the theorem.
\end{prf}

\subsubsection{Noncritical level}
Consider the algebra $\wtld{\Ac}_{c}(\g_N)$ at the noncritical level $c\in\CC$, $c\neq -N/2$.  Introduce the family of variables $u_-=(u,u-h,\ldots ,u-(N-1)h)$ and  
denote by $A^{(N)}$  the action of the  anti-symmetrizer from the group algebra $\CC[\Sym_N]$ on the tensor product space $(\CC^N)^{\ot N}$.
Let $\sss  (u)$ be the Laurent series with coefficients in $\wtld{\Ac}_{c}(\g_N)$,
\beq\label{yryxc1}
  \sss  (u)=\tr_{1,\ldots,N}\ts   A^{(N)}   S^+_{[N]}(u_-)\ts S_{[N]}(u_- +hc/2)^{-1}. 
\eeq
Note that $u_-=u_\nu$ if $\nu$ is the column diagram with $N$ boxes; recall \eqref{variable}. Moreover, we have $A^{(N)}=\Ec_\Uc$ if  $\Uc$ is the unique standard column tableaux with $N$ boxes.  
In parallel with the critical level case, we observe  that the action of the series $\sss  (u)$ on the vacuum module $\Wc_{c}(\g_N)$ is given by
\beq\label{yryxc2}
\sss  (u)=Y_{\Wc_{c}}(\qdet T^+(0),u).
\eeq
We have the following simple construction of the family of central elements in $\wtld{\Ac}_{c}(\g_N)$.

\begin{pro}\label{proorp}
All coefficients of  $ \sss  (u)$ belong to the center of the algebra $\wtld{\Ac}_{c}(\g_N)$. 
\end{pro}

\begin{prf}
A direct calculation which relies on the properties of the anti-symmetrizer  and goes in parallel with the
proof of \cite[Thm. 2.5.3]{M}
shows that there exists a power series $\gamma(u)\in\CC[u^{-1}][[h]]$ such that
we have
\begin{align*}
\sss  (u)=&\,\, \gamma(u)\ts \qdet T^+ (u) \ts \qdet T^+(-u+(N-1)h)  \\
  &\times \left(\qdet T(-u-hc/2+(N-1)h)\right)^{-1}   \left(\qdet T(u+3hc/2)\right)^{-1}.
\end{align*}
However, it is well-known that the coefficients of these quantum determinants  
belong to the center of the double Yangian; see, e.g., \cite[Prop. 2.8]{JKMY}. Thus, we conclude that   all coefficients of    $ \sss  (u)$ belong to the center of the algebra $\wtld{\Ac}_{c}(\g_N)$, as required. 
\end{prf}

\subsection{Invariants of the vacuum module \texorpdfstring{$\Wc_{c}(\g_N)$}{Mc(gN)}}
In this section, we present some simple applications of Theorem \ref{main} and Proposition \ref{proorp}, in particular,  
to
the submodule of {\em invariants} of the vacuum module  $\Wc_{c}(\g_N)$, which we define  as
$$\z(\Wc_{c}(\g_N))=\left\{w\in\Wc_{c}(\g_N)\,:\,  S(u)w=G\ot w  \right\}.$$
Our definition is motivated by the usual notion of a subspace of invariants of the vacuum module for the universal enveloping algebra of   affine Lie algebra; see, e.g., the book by E. Frenkel \cite{F}.
As before, we shall consider the critical and the noncritical level separately. It is worth noting  that for any $c\in\CC$ we have $1\in \z(\Wc_{c}(\g_N))$ as $S(u) 1 = G\ot 1$. 

\subsubsection{Critical level}
Consider the vacuum module $\Wc_{-N/2}(\g_N)$ at the critical level $c=-N/2$.
By applying the   series \eqref{6c1}, whose coefficients belong to the center of the algebra $\wtld{\Ac}_{-N/2}(\g_N)$, on $1\in \Wc_{-N/2}(\g_N)$ we obtain  
$$
  \mmm_{\nu }  (u)=\sss_{\nu } (u)1 \in \Wc_{-N/2}(\g_N)[[u^{\pm 1}]]. 
$$

\begin{kor}\label{Bns}
All coefficients of the series $\mmm_{\nu }  (u)$  belong to   $\z(\Wc_{-N/2}(\g_N))$. 
\end{kor}
\begin{prf}
The corollary is a simple consequence of  
 Theorem \ref{main}.  
We have 
\begin{align*}
S(v)\ts \mmm_{\nu }  (u) 
= S(v)\ts \sss_{\nu } (u)1
= \sss_{\nu }(u)\ts S(v)  1
= \sss_{\nu }(u) (G \ot 1)
= G\ot \ts\sss_{\nu }(u)   1
=G \ot\ts\mmm_{\nu }(u),
\end{align*}
so the coefficients of   $\mmm_{\nu }  (u)$  belong to the submodule of invariants  $\z(\Wc_{-N/2}(\g_N))$. 
\end{prf}

Note that the coefficients of the series $\mmm_{\nu }  (u)$ can be naturally regarded as elements of the $h$-adically completed algebra
$\Y^+ (\g_N)_h$. This leads to   another simple consequence of Theorem \ref{main}, which produces commutative subalgebras in $\Y^+ (\g_N)_h$.

\begin{kor}\label{bnssnb}
The coefficients of all series $\mmm_{\nu } (u)\in \Y^+ (\g_N)_h [[u^{\pm 1}]]$  pairwise commute.
\end{kor}

\begin{prf}
Let $\nu_i$ with $i=1,2$ be   partitions with at most $N$ parts. By Theorem \ref{main} we have
$$
 \mmm_{\nu_i }  (u) \ts \mmm_{\nu_j }  (v)
= \mmm_{\nu_i }  (u) \ts \sss_{\nu_j }  (v)1
=\sss_{\nu_j }  (v)\ts \mmm_{\nu_i }  (u)1
=\sss_{\nu_j }  (v)\ts \sss_{\nu_i }  (u)1.
$$
for all $i,j=1,2$. Since $\sss_{\nu_i }  (u)$ and $\sss_{\nu_j }  (v)$ commute, this implies the corollary.
\end{prf}

Let $\mu$ and $\nu$ be   partitions  which have at most $N$ parts.
Theorem \ref{main}   implies the identity
$ 
\sss_{\nu  }  (u)\ts \sss_{\mu  }  (v) =\sss_{\mu  }  (v)\ts \sss_{\nu }  (u)
$ 
for operators on   $\Wc_{-N/2}(\g_N)$. Hence, by using \eqref{sfh1} we find
\beq\label{sfh2}
Y_{\Wc_{-N/2} }(\TT^{+}_{\nu} (u),z_1)\ts 
Y_{\Wc_{-N/2} }(\TT^{+}_{\mu} (v),z_2) 
= 
Y_{\Wc_{-N/2} }(\TT^{+}_{\mu} (v),z_2) \ts
Y_{\Wc_{-N/2} }(\TT^{+}_{\nu} (u),z_1)
.
\eeq
This commutation relation can be generalized as follows.

\begin{kor}\label{main2}
For any two elements $a$ and $b$ of the center $\z(\Vc_{-N}(\gl_N))$   we have
\beq\label{kom4}
Y_{\Wc_{-N/2} }(a,z_1)\ts
Y_{\Wc_{-N/2} }(b,z_2)\ts
=
Y_{\Wc_{-N/2} }(b,z_2)\ts
Y_{\Wc_{-N/2} }(a,z_1) .
\eeq
\end{kor}

\begin{prf}
By \cite[Thm. 4.9]{JKMY}
and  \eqref{sfh2}, there exists a family of topological generators of the center $\z(\Vc_{-N}(\gl_N))$ such that their images under the quasi module map  \eqref{Y} commute. However, if the commutativity relation  \eqref{kom4} holds for all elements $a$ and $b$ of the family of topological generators of the center, then it holds for all elements of the center as well; see   \cite[Prop. 3.4]{K19}.
Hence
 the identity \eqref{kom4} holds for any $a,b\in \z(\Vc_{-N}(\gl_N))$.
\end{prf}

We conclude this section by the   method for constructing invariants of   $\Wc_{-N/2}(\g_N)$.

\begin{kor}\label{zadnji}
For any $a\in\z(\Vc_{-N}(\gl_N))$ and $w\in \z(\Wc_{-N/2}(\g_N))$ all coefficients of the series $Y_{\Wc_{-N/2} }(a,z)w$ belong to the submodule of invariants  $\z(\Wc_{-N/2}(\g_N))$.
\end{kor}

\begin{prf}
As with the proof of Corollary \ref{main2}, we can combine \cite[Thm. 4.9]{JKMY}, \cite[Prop. 3.4]{K19}
and  \eqref{sfh2} to conclude that for  operators on $\Wc_{-N/2}(\g_N)$  we have
$$
Y_{\Wc_{-N/2} }(a,z)\ts S(u)=S(u)\ts Y_{\Wc_{-N/2} }(a,z)\quad\text{for all}\quad a\in\z(\Vc_{-N}(\gl_N)) .
$$
Therefore, for any $a\in\z(\Vc_{-N}(\gl_N))$ and $w\in \z(\Wc_{-N/2}(\gl_N))$ we have
\begin{align*}
S(u)\ts Y_{\Wc_{-N/2} }(a,z)w =Y_{\Wc_{-N/2} }(a,z)\ts  S(u)w=Y_{\Wc_{-N/2} }(a,z) (G\ot w) =G\ot Y_{\Wc_{-N/2} }(a,z)w.
\end{align*}
Hence the coefficients of $Y_{\Wc_{-N/2} }(a,z)w$ belong to   $\z(\Wc_{-N/2}(\g_N))$, as required.
\end{prf}

\subsubsection{Noncritical level}
Consider the vacuum module $\Wc_{c}(\g_N)$ at the noncritical  level $c\in\CC$, $c\neq -N/2$.
As with the critical level case,   applying the   series \eqref{yryxc1}, whose coefficients belong to the center of  $\wtld{\Ac}_{c}(\g_N)$, on $1\in \Wc_{c}(\g_N)$ we obtain a   power series
$$
  \mmm   (u)=\sss  (u)1 \in \Wc_{c}(\g_N)[[u^{\pm 1}]]. 
$$
 By  arguing as in   the proof of Corollary \ref{Bns}  and using Proposition \ref{proorp} one  obtains
\begin{kor} 
All coefficients of the series $\mmm   (u)$  belong to   $\z(\Wc_{c}(\g_N))$. 
\end{kor}
 The analogues of Corollaries  \ref{bnssnb}, \ref{main2} and \ref{zadnji} for the series $\mmm   (u)$  can be   easily established as well. However, their proofs   now rely on the Proposition \ref{proorp} and the explicit description of the center $\z(\Vc_{c}(\gl_N))$ at the noncritical level $c\neq -N$, as given by \cite[Prop. 4.10]{JKMY}.

\begin{rem}
It is worth to single out the following   identity  for operators on $\Wc_{c}(\g_N)$:
\begin{align*}
Y_{\Wc_{c} }(\qdet T^+(0),z)=
\sdet S^+  (z) \ts \sdet S( z+ hc/2)^{-1}. 
\end{align*}
It connects the quantum determinant \eqref{kvdet} and the {\em Sklyanin determinants}, 
\begin{align*}
\sdet S^+ (u)=\tr_{1,\ldots,N}\ts   A^{(N)}   S^+_{[N]}(u_-)\fand 
\sdet S  (u)=\tr_{1,\ldots,N}\ts   A^{(N)}   S _{[N]}(u_-).
\end{align*}
The Sklyanin determinant   $\sdet S^+ (u)$ exhibits   similar properties as its  more intensively studied Yangian counterpart $\sdet S (u)$; cf. \cite{O,MNO,M}. 
\end{rem}

\section*{Acknowledgement}
This work has been supported in part by Croatian Science Foundation under the project UIP-2019-04-8488.

\linespread{1.0}

\end{document}